\input amstex
\documentstyle{amsppt}
\nopagenumbers \accentedsymbol\tx{\tilde x}
\TagsOnRight
\accentedsymbol\ty{\tilde y}

\pagewidth{13cm} 
\pageheight{21cm} 
\rightheadtext{Point-invariant classes \dots} 
\topmatter
\title
Point-invariant classes of the third order ordinary 
differential equations.
\endtitle
\author
V.~V.~Dmitrieva
\endauthor
\address
Department of Mathematics, Bashkir State University, Frunze str\.~32,
450074 Ufa, Russia.
\endaddress
\email
DmitrievaVV\@ic.bashedu.ru
\endemail
\thanks
Author is grateful to Prof\.~R.~A.~Sharipov for stating of the problem and for
useful discussions. This research was supported by RFBR, grant no\. 00-01-00068,
coordinator Prof.~Y.~T.~Sultanaev.
\endthanks

\abstract Point transformations of the 3-rd order ordinary differential
equations are considered. Special
classes of ordinary differential equations
that are invariant under the general point transformations are
constructed.
\endabstract

\endtopmatter
\loadbold
\document

\head
1. Introduction.
\endhead

Let us consider a class of ordinary differential equations of the third order
resolved with respect to the leading derivative $$ y^{\prime\prime\prime}=
f(x,\,y,\,y^\prime,\,y^{\prime\prime}).\hskip -2em\tag{1.1} $$ Under the general
point transformation
$$\hskip -2em \cases
\tilde x=\tilde x(x,y),\\
\tilde y=\tilde y(x,y).
\endcases
\tag{1.2} $$
some equation from the class \thetag{1.1} is transformed into the
following one
$$ \tilde y^{\prime\prime\prime}=g(\tilde x,\,\tilde y,\,\tilde
y^\prime,\,\tilde y^{\prime\prime}).\hskip -2em \tag1.3
$$

Let us choose the special subclass of ODE's \thetag{1.1} such that the function $f$
from the right hand side of \thetag{1.1} is in the certain preassigned set of
functions $\Cal F$ ($f\in {\Cal F}$).

\proclaim {Definition 1}{If $\,\,\forall f\in {\Cal F}$ the function $g$ from the
right hand side of transformed equation \thetag{1.3} is also in the set $\Cal F$
($g\in {\Cal F}$) then this subclass \thetag{1.1} is called a point-invariant
subclass.
\endproclaim

As usual one consider the class $\Cal F$ consisting of polynomials or rational
functions in derivatives $y'$ and $y''$.

In XIX-th century the following ODE class $$ \hskip -2em
y^{\prime\prime}=P(x,y)+3\,Q(x,y)\,y^\prime+3\,R(x,y)\,(y^\prime)^2+
S(x,y)\,(y^\prime)^3.\hskip -2em \tag{1.4} $$ was investigated by R.~Liouville
\cite{1}, S.~Lie \cite{2}-\cite{3} and A.~Tresse \cite{4}-\cite{5}. E.~Cartan
considered \thetag{1.4} as the equation of geodesic line in the projectively
connected space \cite{6} -- \cite{9}. The class of equations \thetag{1.4} is closed
under the point transformations \thetag{1.2}. This fact was mentioned in
\cite{1}--\cite{24}.

In \cite{25} R.~A.~Sharipov and O.N.~Mikhailov have investigated some other
point-invariant class of second order equations: $$ y^{\prime\prime}= \frac
{P(x,y)+4\,Q(x,y)\,y^\prime+6\,R(x,y)\,(y^\prime)^2+
4S(x,y)\,(y^\prime)^3+L(x,y)(y^\prime)^4} {Y(x,y)-X(x,y)y^\prime}.\hskip -2em
$$ They called it {\it the point-expansion} for the class of equations
\thetag{1.4}.

Another class was mentioned L.A.~Bordag in \cite{17}: $$ (y'')^2=P_5(y';x,y). $$
Here $P_5=P_5(y';x,y)$ is a polynomial of the fifth order in $y'$.

\head
2. Point transformations.
\endhead

We can interpret the change of variables \thetag{1.2} as the change of some
curvilinear coordinates on the plane for another ones. Let's denote by T and S the
direct and inverse Jacoby matrices for this change of variables: $$ \xalignat 2 &S =
\Vmatrix {x_{\sssize 1. 0}} & {x_{\sssize 0. 1}} \\ \vspace{1ex} {y_{\sssize 1. 0}}
& {y_{\sssize 0. 1}}
\endVmatrix,
&& T =  \Vmatrix {\tilde x_{\sssize 1.0}} & {\tilde x_{\sssize 0.
1}} \\ \vspace{1ex} {\tilde y_{\sssize 1. 0}} & {\tilde y_{\sssize
0.1}}
\endVmatrix.
\tag{2.1}
\endxalignat
$$ The point transformations \thetag{1.2} are supposed to be non-degenerate, hence
$\det S\ne 0$. Here and below the notation $\Phi_{\sssize i.j}$ denotes the partial
derivative: $$\Phi_{\sssize i.j}=\frac {\partial^{i+j}\Phi }{\partial x^i\partial
y^j}.$$ In particular $$x_{\sssize 1. 0}=\frac {\partial x(\tilde x, \tilde y)}
{\partial\tilde x}.$$

The derivative $y'$ obeys the following transformation rule under the change of
coordinates \thetag{2.1}: $$ {y^{\prime }} =  \frac {{y_{\sssize 1. 0}} +
{y_{\sssize 0.1}}\,{\tilde y^{\prime }}}{{x_{\sssize 1. 0}} + {x_{\sssize 0. 1}}\,
{\tilde y^{\prime }} }.\tag{2.2} $$

The corresponding rule for the $y''$ is more complicated:
$$ \aligned
{y^{\prime\prime }} & =   \frac { {y_{\sssize 2. 0}} +
2\,{y_{\sssize 1. 1}}\,{\tilde y^{\prime }} + {y_{ \sssize
0.2}}\,{\tilde y^{\prime2 }} + {y_{\sssize 0. 1}}\,{\tilde y
^{\prime\prime }}}{({x_{\sssize 1. 0}} + {x_{\sssize 0.1}}\,
{\tilde y ^{\prime }})^2 } -\\
 & -  \frac {({y_{\sssize 1. 0}} + {y_{\sssize 0. 1}}
\, {\tilde y^{\prime }})\,({x_{\sssize 2. 0}} + 2\,{x_{\sssize 1.
1}}\, {\tilde y^{\prime }} + {x_{\sssize 0. 2}}\, \tilde y^{\prime
2}
 + {x_{\sssize 0. 1}}\, {\tilde y^{\prime\prime }} )} {({x_{
\sssize 1. 0}}
 + {x_{\sssize 0. 1}}\, \tilde y^{\prime })^3 }
\endaligned\tag{2.3}
$$

Appropriate formula for the $y'''$: $$ \aligned y^{\prime\prime\prime}= & \frac
{a_1\tilde y^{\prime\prime\prime}+ a_2\tilde y^{\prime\prime 2}+a_3\tilde
y^{\prime\prime} {\tilde y^{\prime 2}}+a_4{\tilde y^{\prime\prime}}\,{\tilde
y^\prime}+ a_5{\tilde y^{\prime\prime}}+} {({x_{\sssize 1. 0}} + {x_{\sssize 0.
1}}\,{\tilde y^{\prime}})^5} \\ & \frac {+a_6{\tilde y^{\prime 5}}+ a_7{\tilde
y^{\prime 4}}+a_8{\tilde y^{\prime 3}}+ a_9{\tilde y^{\prime 2}}+a_{10}{\tilde
y^\prime}+a_{11}} {({x_{\sssize 1. 0}} + {x_{\sssize 0. 1}}\,{\tilde y^\prime})^5}.
\endaligned\tag {2.4}
$$

All coefficients  $a_2,\dots, a_{11}$ in \thetag {2.4} are some certain functions
in  $x$, $y$ and their derivatives with respect to $\tilde x$, $\tilde y$, except
for the coefficient $a_1$ which contains the dependence on $\tilde y'$. Let's write
the explicit formulae for $a_1,\,\dots,\,a_{11}$.

The coefficient of $\tilde y^{\prime\prime\prime}$: $$ a_1=  - ({x_{\sssize 1.
0}}   +   {x_{\sssize 0. 1}}\,{\tilde y^\prime})\,({y_{\sssize 1. 0}}{x_{\sssize 0.
1}} - {x_{\sssize 1. 0}}\,{y_{\sssize 0. 1}}) \tag {2.5} $$

The coefficient of $\tilde y^{\prime\prime 2}$: $$ a_2= 3{x_{\sssize 0.
1}}\,({y_{\sssize 1. 0}}\,{x_ {\sssize 0. 1}} - {x _{\sssize 1. 0}}\,{y_{\sssize 0.
1}})\tag {2.6} $$

The coefficient of $\tilde y^{\prime\prime}\tilde y^{\prime 2}$:  $$ \aligned a_3=&
- 6\,{y_{\sssize 0. 1}}\,{x_{\sssize 1. 0}}\,{x_{\sssize 0 .2}} - 3\,{y_{\sssize 1.
1}}\,{x_{\sssize 0. 1}}^{2} + 3\,{y_ {\sssize 1. 0}}\,{x _{\sssize 0.
2}}\,{x_{\sssize 0. 1}} + \\ & + 3\,{y_{\sssize 0. 2}}\,{x_{\sssize 1.
0}}\,{x_{\sssize 0. 1}} + 3\,{y_{\sssize 0. 1}}\,{x_{\sssize 1. 1}}\,{x_{\sssize
0.1}}
\endaligned\tag{2.7}
$$

The coefficient of $\tilde y^{\prime\prime}\tilde y^{\prime }$: $$ \aligned a_4= &
9\,{y_{\sssize 1. 0}}\,{x_{\sssize 0. 1}}\,{x_{\sssize 1. 1}} - 3\,{y_{\sssize 2.
0}}\,{x_{\sssize 0. 1}}^{2} + 3\,{y_{ \sssize 0. 1}}\,{x_{ \sssize 2.
0}}\,{x_{\sssize 0. 1}} +  \\ & 3\,{y_{\sssize 0. 2}}\,{x_{\sssize 1. 0}}^{2} -3\,{
y_{\sssize 1. 0}}\,{x_{\sssize 0.2}}\,{x_{\sssize 1.0}} -9\,{y_{\sssize
0.1}}\,{x_{\sssize 1. 1}}\,{x_{\sssize 1.0}}
\endaligned\tag{2.8}
$$

The coefficient of $\tilde y^{\prime\prime}$: $$ \aligned a_5= & - 3\,{y_{\sssize 1.
0}}\,{x_{\sssize 1. 1}}\,{x_{\sssize 1.0}} + 6\,{y_{1,0}}\,{x_{\sssize 2.
0}}\,{x_{\sssize 0. 1}}-
  3\,{y_{\sssize 0. 1}}\,{x_{\sssize 1.0}}
\,{x_{\sssize 2. 0}} + \\ & +3\,{y_{\sssize 1. 1}}\,{x_{\sssize 1.
0}}^{2}  -  3\,{y_{\sssize 2.0}}\,{x_{\sssize 1. 0}}\,{x_{\sssize
0.1}}
\endaligned\tag{2.9}
$$

The coefficient of $\tilde y^{\prime 5}$: $$ \aligned a_6= & 3\,{y_{\sssize 0.
1}}\,{x_{\sssize 0. 2}}^{2} - 3\,{y_{\sssize 0.2}}\,{x _{\sssize 0.
1}}\,{x_{\sssize 0. 2}} + {y_{\sssize 0.3}}\,{x_{\sssize 0. 1}}^{2} - \\
&-{y_{\sssize 0. 1}}\,{x_{\sssize 0. 3}}\,{x_{\sssize 0.1}}
\endaligned\tag{2.10}
$$

The coefficient of $\tilde y^{\prime 4}$: $$ \aligned a_7 = &  - 3\,{y_{\sssize 0
.1}}\,{x_{\sssize 1. 2}}\,{x_{\sssize 0.1}}
 - 3\,{y_{\sssize 0. 2}}\,{x_{\sssize 1. 0}}\,{x_{\sssize 0. 2}} - {y_{
\sssize 0. 1}}\,{x_{\sssize 0. 3}}\,{x_{\sssize 1. 0}} - \\
&-6\,{y_{\sssize 1. 1}}\,{x_{\sssize 0. 2}}\,{x_{\sssize 0. 1}} +
2\,{y_{\sssize 0. 3}}\,{x_{\sssize 1. 0}}\, {x_{\sssize 0. 1}}
 + 3\,{y_{\sssize 1. 0}}\,{x_{\sssize 0. 2}}^{2} -\\
& - {y_{\sssize 1. 0}}\,{x_{\sssize 0. 3}}\,{x_{\sssize 0. 1}} -
6\,{y_{ \sssize 0.2}}\,{x_{\sssize 1. 1}}\,{x_{\sssize 0. 1}} +
12\,{y_{\sssize 0. 1}}\,{x_{\sssize 1. 1}}\,{x_{\sssize 0. 2}}  +
3\,{y_{\sssize 1. 2}}\,{x_{\sssize 0.1}}^{2}
\endaligned\tag{2.11}
$$

The coefficient of $\tilde y^{\prime 3}$: $$ \aligned a_8=& - 3\,{y_{\sssize 0.
1}}\,{x_{\sssize 2.1}}\,{x_{\sssize 0. 1}} - 6\,{y_{\sssize 0. 2}}\,{x_{\sssize 1.
1}}\,{x_{\sssize 1.0}} -
{y_{\sssize 1. 0}}\,{x_{\sssize 0. 3}}\,{ x_{\sssize 1. 0}} - \\
&-3\,{y_{\sssize 0. 1}}\,{x_{\sssize 1. 2}}\,{x_{\sssize 1.0}} +
12\,{ y_{\sssize 0. 1}}\,{x_{\sssize 1. 1}}^{2} - 3\,{ y_{\sssize
0.2}}\,{x_{\sssize 2.0}}\,{x_{\sssize 0.1}}-\\
 &   - 3\,{y_{\sssize 2. 0}}\,{x_{\sssize 0
.2}}\,{x_{\sssize 0. 1}}
 - 6\,{y_{\sssize 1.1}}\,{x_{\sssize 0. 2}}\,{x_{\sssize 1.0
}} + 6\,{y_{\sssize 0. 1}}\,{x_{\sssize 2. 0}}\,{x_{\sssize 0.
2}}-\\ &  - 12\,{y_{\sssize 1.1}}\,{x_{\sssize 0.1}}\,{x_{\sssize
1. 1}} + 6\, {y_{\sssize 1.2}}\,{x_{\sssize 1. 0}}\,{x_{\sssize
0.1}} + {y_{\sssize 0. 3}}\,{x_{\sssize 1. 0}}^{2} + \\ & +
3\,{y_{\sssize 2. 1}}\,{x_{\sssize 0. 1}}^{2}  -  3\,{y_{\sssize
1. 0}}\, {x_{\sssize 1. 2}}\,{x_{\sssize 0. 1}}
  + 12\,{y_{\sssize 1. 0}}\,{x_{\sssize 1. 1}}\,{x_{\sssize 0. 2}}
\endaligned\tag{2.12}
$$

The coefficient of $\tilde y^{\prime 2}$: $$ \aligned a_9 =& - 3\,{y_{\sssize 1.
0}}\,{x_{\sssize 1. 2}}\,{x_{\sssize 1.0}}
 + 12\,{y_{\sssize 0. 1}}\,{x_{\sssize 2. 0}}\,{x_{\sssize 1. 1}} +
3\,{y_{\sssize 1. 2}} \,{x_{\sssize 1. 0}}^2  +\\
 & + {y_{\sssize 3. 0}}\,{x_{\sssize 0.1}}^{2} - 6\,{y_{\sssize 2. 0}}\,
{x_{\sssize 1. 1}}\,{x_{\sssize 0.1}} - 3\,{y_{\sssize
0.1}}\,{x_{\sssize 2. 1}}\,{x_{\sssize 1 .0}} - \\ &
-6\,{y_{\sssize 1. 1}}\,{x_{\sssize 2.0}}\,{x_{\sssize 0.1}} -
3\,{y_ {\sssize 1 .0}}\,{x_{\sssize 2. 1}}\,{x_{\sssize 0. 1}} -
{y_{\sssize 0.1}}\,{x_{\sssize 3. 0}}\,{x_{\sssize 0. 1}} -\\
 &  -3\,{y_{\sssize
0.2}}\,{x_{\sssize 2. 0}}\,{x_{\sssize 1. 0}} + 6\,{y_{\sssize
1.0}}\, {x_{\sssize 2. 0}}\,{x_{\sssize 0. 2}} + 6\,{y_{\sssize 2.
1}}\,{x_{\sssize 1. 0}}\, {x_{\sssize 0. 1}}-\\ & -
12\,{y_{\sssize 1. 1}}\,{x_{\sssize 1. 0}}\,{x_{\sssize 1. 1}}
 + 12\,{y_{\sssize 1. 0}}\,{x_{\sssize 1. 1}}^{2}  -  3\,{y_{\sssize
2. 0 }}\,{x_{\sssize 0. 2}}\,{x_{\sssize 1. 0}}.
\endaligned\tag{2.13}
$$

The coefficient of $\tilde y^{\prime}$: $$ \aligned a_{10}=& - 3\,{y_{\sssize 2
.0}}\,{x_{\sssize 0. 1}}\,{x_{\sssize 2. 0}} +
12\,{y_{\sssize 1. 0}} {x_{\sssize 2. 0}}{x_{\sssize 1. 1}} -\\
 &  - 6\,{y_{\sssize 2. 0}}\,{x_{\sssize 1.1}}\,{x_{\sssize 1. 0}} -
6\,{ y_{\sssize 1.1}}\,{x_{\sssize 2. 0}}\,{x_{\sssize 1. 0}} -
3\,{y_{\sssize 1. 0}}\,{x_{\sssize 2 . 1}}\,{x_{\sssize 1.0}} - \\
& -{y_{\sssize 0.1}}\,{x_{\sssize 3.0}}\,{x_{\sssize 1. 0}}
 + 2\,{y_{\sssize 3.0}}\,{x_{\sssize 1. 0}}\,{x_{\sssize 0. 1}} +\\
&  + 3\,{y_{\sssize 0.1}}\,{x_{\sssize 2. 0}}^{2}  + 3\,{y_{\sssize 2.1}
}\,{x_{\sssize 1. 0}}^{2} - {y_{\sssize 1.0}}\,{x_{\sssize 3.  0}}\, {x_{\sssize
0.  1}}.
\endaligned\tag{2.14}
$$

The last term: $$ \aligned a_{11}=& - {y_{\sssize 1. 0}}\,{x_{\sssize 3.
0}}\,{x_{\sssize 1.0}} + {y_{\sssize 3.
0}}\,{x_{\sssize 1. 0}}^{2} +\\
 &  + 3\,{y_{\sssize 1. 0}}\,{x_{\sssize 2.0}}^{2} -  3\,{y_{\sssize
2.0} }\,{x_{\sssize 1.0}}\,{x_{\sssize 2. 0}}.
\endaligned\tag{2.15}
$$

Formula \thetag{2.3} and the series of formulae \thetag{2.4} -- \thetag{2.15}
determine the  transformation rule for $y'''$ under the change of coordinates
\thetag {2.1}.

\head 3. Point-invariant classes of equations of the form \thetag{1.1}.
\endhead

Below we shall seek the point-invariant classes of the  equations of the form
\thetag{1.1} such that the functions $f$ in their right hand sides are rational in
$y'$ and polynomial in $y''$. At first we apply the transformation rules
\thetag{2.4}--\thetag{2.15} to the simplest equation $$ y'''=0. \tag{3.1} $$

The shape of transformed equation \thetag{3.1} $$ \aligned \tilde
y^{\prime\prime\prime}= & \frac {3x_{\sssize 0.1}\tilde y^{\prime\prime 2}}
{x_{\sssize 1.0}+x_{\sssize 0.1}\tilde y'} +\frac {a_3\tilde
y^{\prime\prime}{\tilde y^{\prime 2}}+ a_4{\tilde y^{\prime\prime}}\,{\tilde
y^\prime}+a_5{\tilde y^{\prime\prime}}+} {({x_{\sssize 1. 0}} + {x_{\sssize
0.1}}\,{\tilde y^{\prime}}) (y_{\sssize 1.0}x_{\sssize 0.1}-x_{\sssize
1.0}y_{\sssize 0.1})} \\ & \frac {+a_6{\tilde y^{\prime 5}}+ a_7{\tilde y^{\prime
4}}+a_8{\tilde y^{\prime 3}}+ a_9{\tilde y^{\prime 2}}+a_{10}{\tilde
y^\prime}+a_{11}} {({x_{\sssize 1. 0}} + {x_{\sssize 0.1}}\,{\tilde y^\prime})
(y_{\sssize 1.0}x_{\sssize 0.1}-x_{\sssize 1.0}y_{\sssize 0.1})}
\endaligned
\tag {3.2} $$ allows us to assume that the point-invariant class of the equations
\thetag{1.1} should have the following form $$ \aligned y^{\prime\prime\prime}=
&\frac { B(x,y) y^{\prime\prime 2}+P(x,y) y^{\prime\prime} { y^{\prime 2}}+Q(x,y){
y^{\prime\prime}}\,{ y^\prime}+ R(x,y){ y^{\prime\prime}}+S(x,y){y^{\prime 5}}+}
{Y(x,y)-X(x,y){ y^{\prime}}} \\ & \frac {+L(x,y){ y^{\prime 4}}+K(x,y){ y^{\prime
3}}+ M(x,y){ y^{\prime 2}}+N(x,y){ y^\prime}+T(x,y)} {Y(x,y)-X(x,y){ y^\prime}}
\endaligned.\tag{3.3}
$$

Let's transform an arbitrary  equation  \thetag{3.3} using formulae
\thetag{2.2}--\thetag{2.15}. In general the transformed equation do not belong to
the class \thetag{3.3}.  The coefficient $B/(Y-Xy')$  of  $y^{\prime\prime 2}$ in
the right hand side of \thetag{3.3} as a  result of coordinate transformations
\thetag{1.2} generates a new term with $y'y^{\prime\prime 2}$. $$ \tilde y'''=
\frac {(\tilde B_1+\tilde B_2\tilde y')\tilde y^{\prime\prime 2}} {(\tilde Y-\tilde
X\tilde y')(x_{\sssize 1.0}+x_{\sssize 0.1}\tilde y')} +\dots.\tag{3.4} $$

Let us introduce a new function $f(z)$ of the variable $z=y'$ associated with the
equation \thetag{3.4}: 
$$ f(z)=\frac {\tilde B_1+\tilde B_2 z} {x_{\sssize
1.0}+x_{\sssize 0.1}z}, \tag{3.5} $$ 
where 
$$ \aligned &\tilde B_1=3x_{\sssize
0.1}(Yx_{\sssize 0.1}-Xy_{\sssize 0.1})\\ 
&\tilde B_2=3x_{\sssize 0.1}(Yx_{\sssize
1.0}-Xy_{\sssize 1.0})+ B(x_{\sssize 1.0}
y_{\sssize 0.1}-x_{\sssize 0.1}y_{\sssize 1.0}).
\endaligned
\tag{3.6}
$$

The rational function $f(z)$ in \thetag{3.5} has a first order pole at the point
$z_0=-x_{\sssize 1.0}/x_{\sssize 0.1}$. Let's calculate the residue of this
function at the point  $z_0$ and denote it by $\Omega$: $$ \Omega=
\operatornamewithlimits{Res} _{z=z_0} f(z).\tag{3.7} $$

It is easy to get the explicit formula for $\Omega$ using the expression
\thetag{3.6}: $$ \Omega=(B+3X)\det S. \tag{3.8} $$

The condition $\Omega=0$ is the  {\it necessary condition} for the class of
equations \thetag{3.2} to be point-invariant. Let's remember that the point
transformations \thetag{1.2} are non-degenerate, hence the condition \thetag{3.8} is
equivalent  to additional relation between functions $B$ and $X$: $$ B=-3X.\tag{3.9}
$$

The direct calculations show that the condition  \thetag{3.9} is also {\it
sufficient condition } for the class of equations \thetag{3.2} to be closed
 under the transformations \thetag{1.2}. This completes the proof of Theorem 1
 below.

\proclaim{Theorem 1} {The class of equations
$$ \aligned
y^{\prime\prime\prime}= & \frac { -3X(x,y) y^{\prime\prime
2}+P(x,y) y^{\prime\prime} { y^{\prime 2}}+Q(x,y){
y^{\prime\prime}}\,{ y^\prime}+ R(x,y){
y^{\prime\prime}}+S(x,y){y^{\prime 5}}+} {Y(x,y)-X(x,y){
y^{\prime}}} \\ & \frac {+L(x,y){ y^{\prime 4}}+K(x,y){ y^{\prime
3}}+ M(x,y){ y^{\prime 2}}+N(x,y){ y^\prime}+T(x,y)}
{Y(x,y)-X(x,y){ y^\prime}}
\endaligned\tag{3.10}
$$ is invariant  under point transformations \thetag{1.2}.}
\endproclaim







\Refs
\ref\no 1\by R.~Liouville \paper Sur les invariants de
certaines equations differentielles et sur leurs applications
\jour J. de L'Ecole Polytechnique \vol 59 \yr 1889 \pages 7--76
\endref

\ref\no 2\by S.Lie \book Vorlesungen \"uber continuierliche
Gruppen \publ Teubner Verlag \publaddr Leipzig \yr 1893
\endref

\ref\no 3\by S.Lie \book Theorie der Transformationsgruppen III
\publ Teubner Verlag \publaddr Leipzig \yr 1930
\endref

\ref\no 4\by A.~Tresse \paper Sur les invariants differenties des
groupes continus de transformations \jour Acta Math. \vol 18
\pages 1--88 \yr 1894
\endref

\ref\no 5 \by A.~Tresse\paper Determination des Invariants
ponctuels de l"Equation differentielle ordinaire de second ordre:
$y''=w(x,y,y')$. \inbook Preisschriften der f\"rstlichen
Jablonowski'schen Gesellschaft XXXII \publ S.Hirzel\publaddr
Leipzig \yr 1896
\endref

\ref\no 6\by E.~Cartan\paper Sur les varietes a connection
projective\jour Bulletin de Soc. Math. de France \vol 52 \yr
1924\pages 205--241
\endref

\ref\no 7\by E.~Cartan\paper Sur les varietes a connexion affine
et la theorie de la relativite generalisee\jour Ann. de l'Ecole
Normale,\vol 40\pages 325--412\yr 1923\moreref\vol 41\yr 1924
\pages 1--25\moreref\yr 1925\vol 42\pages 17--88
\endref

\ref\no 8\by E.~Cartan\paper Sur les espaces a connexion
conforme\jour Ann. Soc. Math. Pologne \vol 2\yr 1923 \pages
171--221
\endref

\ref\no 9\by E.~Cartan\book Spaces of affine, projective and
conformal connection  \publ Platon \publaddr Volgograd\yr 1997
\endref

 \ref \no 10 \by G.~Thomsen \paper {\"U}ber die topologischen
Invarianten der Differentialgleichung $
y''=f(x,y){y'}^3+g(x,y){y'}^2+h(x,y) y' +k(x,y) $ \jour
Abhandlungen aus dem mathematischen Seminar der Hamburgischen
Universit\"at, \vol 7 \yr 1930 \pages 301--328
\endref

\ref\no 11\by C.Grissom, G.Thompson and G.Wilkens \paper
Linearisation of Second Order Ordinary Differential Equations via
Cartan's Equivalence Method \jour Diff. Equations \vol 77\pages
1-15 \yr 1989
\endref

\ref \no 12 \by N.~Kh.~Ibragimov \book Elementary Lie group analisys and ordinary
differential equations, {\rm Wiley series in mathematical methods in practice;
Vol.4} \publ  John Wiley \& Sons Ltd \publaddr England \yr 1999
\endref

\ref\no 13\by V.S.~Dryuma\book Geometrical theory of nonlinear
dynamical system \publ Preprint of Math. Inst. of Moldova
\publaddr Kishinev\yr 1986
\endref

\ref\no 14\by V.S.~Dryuma\paper On the theory of submanifolds of
projective spaces given by the differential equations \inbook
Sbornik statey\publ Math. Inst. of Moldova \publaddr Kishinev\yr
1989\pages 75--87
\endref

\ref\no 15\by V.S.~Dryuma\paper Geometrical properties of
multidimensional nonlinear differential equations and phase space
of dynamical systems with Finslerian metric \jour Theor. and Math.
Phys.,\vol 99\issue 2\yr 1994 \pages 241-249
\endref


\ref\no 16\by L.A.~Bordag and V.S.~Dryuma \paper Investigation of
dynamical systems using tools of the theory of invariants and
projective geometry \inbook NTZ-Preprnt 24/95 \publaddr Leipzig
\yr 1995 \moreref \jour J. of Applied Mathematics (ZAMP) in appear;
\moreref \inbook Electronic archive at LANL (1997), solv-int
\#9705006\pages 1--18
\endref

\ref\no 17\by L.A.~Bordag \paper Symmetries of the Painleve
equations and the connection with projective differential geometry
\inbook VIIth EWM Meeting Proceedings \publ Hindawi Publishing
Corporation \publaddr Trieste, Italy \yr 1997 \pages 145--159
\endref

\ref\no 18\by M.V.~Babich and L.A.~Bordag \paper Projective
Differential Geometrical Structure ot the Painleve Equations\jour
J. of Diff.Equations \vol 157\issue 2, September\yr 1999\pages 452--485
\endref

\ref\no 19\by Yu.R.~Romanovsky\paper Calculation of local
symmetries of second order ordinary differential equations by
means of Cartan's method of equivalence\jour Manuscript \pages
1--20
\endref

\ref\no 20\by S.Bacso and M.Matsumoto\paper On Finsler spaces of Douglas type. A
generalization of the notion of Berwald space \jour Publ. Math. Debrecen \vol
51\issue 3-4\yr 1997 \pages 385--406
\endref

\ref\no 21\by V.N.Gusiatnikova and V.A.Yumaguzhin \paper Point
transformations and linearizability of second-order ordinary
differential equations \jour Mathem. Zametki,\vol 49\issue 1 \yr
1991 \pages 146-148 (in Russian); English transl. in  Soviet Math.
Zametki
\endref


\ref\no 22\by Dmitrieva~V\.~V\., Sharipov~R\.~A\.\paper On the
point transformations for the second order differential
equations\inbook Electronic archive at LANL (1997), solv-int
\#9703003\pages 1--14
\endref

\ref\no 23\by Sharipov~R\.~A\.\paper On the point transformations
for the equation $y''=P+3\,Q\,y'+3\,R\,{y'}^2+S\,{y'}^3$\inbook
Electronic archive at LANL (1997), solv-int \#9706003\pages 1--35
\endref

\ref\no 24\by Sharipov~R.~A.\paper Effective procedure of point
classification for the equations $y''=P+3\,Q\,y'+3\,R\,{y'}^2
+S\,{y'}^3$\inbook Electronic archive at LANL (1998), Math\.DG
\#9802027\pages 1--35
\endref

\ref\no 25\by Mikhailov~O\.~N\., Sharipov~R\.~A\.\paper On the
point expansion for certain class of differential equations of
second order \inbook Electronic archive at LANL (1997), solv-int
\#9712001\pages 1--8
\endref
\endRefs
\enddocument
\end